\documentclass[fleqn]{mat01}
\usepackage{times,mathtimy,amssymb,latexsym}
\begin{document}

\setcounter{page}{55}
\firstpage{65}

\font\xx=msam5 at 10pt
\def\ab{\mbox{\xx{\char'03}}}

\newtheorem{defini}[defin]{\rm DEFINITION}
\def\rem{\trivlist\item[\hskip\labelsep{\it Remark}.]}

\newtheorem{propo}[defin]{\rm PROPOSITION}
\newtheorem{lemm}[defin]{Lemma}
\newtheorem{cor}[defin]{\rm COROLLARY}
\newtheorem{theo}[defin]{\bf Theorem}

\font\ss=tibi at 13.5pt
\font\sa=tibi at 10.4pt

\title{Order units in a $\hbox{\ss C}^{*}$-algebra}

\markboth{Anil K Karn}{Order units in a  $\hbox{\ss C}^{*}$-algebra}

\author{ANIL K KARN}

\address{Department of Mathematics, Deen Dayal Upadhyaya College,
New Delhi 110 015, India\\
\noindent E-mail:~anilkarn@vsnl.net}

\volume{113}

\mon{February}

\parts{1}

\Date{MS received 28 September 2001}

\begin{abstract}
Order unit property of a positive element in a $C^{*}$-algebra is
defined. It is proved that precisely projections satisfy this order
theoretic property. This way, unital hereditary $C^{*}$-subalgebras of a
$C^{*}$-algebra are characterized.
\end{abstract}


\maketitle

\section{Introduction}

Matrix ordered spaces were introduced and studied by Choi and
Effros~\cite{CM}. Matricially normed spaces were introduced and studied
by Ruan~\cite{RZ}. Matrix order unit spaces, studied by Choi and Effros,
are matrix ordered spaces with a matrix norm induced by an order unit.
Their non-unital counterpart were studied by us (see [2--5]) 
and Schreiner \cite{SW}. We call such spaces $L^{\infty}$-
{\it matricially Riesz normed spaces} (Schreiner call them
{\it matrix regular operator spaces}.) We know that $C^{*}$-algebras
are $L^{\infty}$-matricially Riesz normed spaces.

In this paper we consider the following problem: Which element of an
$L^{\infty}$-matricially Riesz normed space looks like order unit and also determines the matrix
norm as the order unit norm in its possible range? (We know that if
two elements determine the (matrix) norm as order unit norms, then the
two elements must be identical.) In this paper, we have been able to
characterize these elements in $C^{*}$-algebras. These are precisely the
projections. This is an order theoretic characterization of projections
in a $C^{*}$-algebra and seems to be an important tool in a possible order theoretic
characterization of non-commutative $C^{*}$-algebras. (Note that the set
of all projections in $B(H)$ have a nice lattice structure.) As a
consequence of the above characterization, we have been able to
characterize unital hereditary $C^{*}$-subalgebras of a $C^{*}$-algebra.

Now we recollect some definitions for the sake of completeness.

\setcounter{defin}{0}
\begin{defini}\label{def1.1}$\left.\right.$\vspace{.5pc}

\noindent {\rm An $L^{\infty}$-{\it matricially Riesz normed space} ($L^{\infty}\hbox{-}mRn$ {\it space})
$V$ consists of a matrix order structure $\{M_{n}(V)^{+}\}$ and a
matrix norm structure $\{\left\|\right\|_{n}\}$ which are
related in the following manner:

\begin{enumerate}
\leftskip .5pc
\renewcommand{\labelenumi}{(\Roman{enumi})}
\item For every $n\in N$ and $v\in M(V)^{+}$
\begin{equation*}
\left\| v\right\|_{n} = \left\lbrace\max\{\left\|u_{1}\right\|_{n}, \left\|u_{2}\right\|_{n}\}: \left(\begin{array}{cc}
u_{1} &v\\
v^{*} &u_{2}
\end{array}\right)\in M_{2n}(V)^{+}\right\rbrace
\end{equation*}
\pagebreak

\item $\left\| v\right\|_{n}\leq \left\| \alpha\right\|\left\|
v\right\|_{n}\left\|\beta \right\|$

\item $M_{n}(V)^{+}$ is $\left\|\right\|_{n}$-closed or all $n$.
\end{enumerate}}
\end{defini}\vspace{-.5pc}

This space is denoted by $(V, \{\left\|\right\|_{n}\},
\{M_{n}(V)^{+}\})$.

Next, we recall the notion of an order unit. Let $(V, \{M_{n}(V)^{+}\})$
be a matrix ordered space. Then $e\in V^{+}$ is called an {\it order
unit} for $V$, if for each $v\in V$ there is $k>0$ such that $\left(\begin{array}{cc}
ke &v\\
v^{*} &ke
\end{array}\right) \in M_{2}(V)^{+}$. In this case, $e^{n} =
e\oplus\cdots \oplus e\in M_{n}(V)^{+}$ for all $n$.
If, in addition, $V^{+}$ is proper and if $M_{n}(V)^{+}$ is Archimedean
for all $n$, then we can define for each $n$, the order unit norm $\|
\cdot \|$ on $M_{n}(V)$ given by
\begin{equation*}
\|v_{n}\| =  {\rm inf}\left\lbrace \lambda > 0: \left\lbrack
\begin{array}{cc}
\lambda e^{n} &v\\
v^{*} &\lambda e^{n}
\end{array}
\right\rbrack  \in M_{2n} (V)^{+}\right\rbrace.
\end{equation*}

In this case $(V, \{\|\cdot\|_{n}\}, \{M_{n} (V)^{+}\})$ becomes an
$L^{\infty}\hbox{-}mRn$ space which is called a matrix order unit space and is
denoted by $(V, e)$.

Now, let $W$ be a self-adjoint subspace of a matrix ordered space $(V,
\{M_{n}(V)^{+}\})$. For each $n\in N$, if we put
\begin{equation*}
M_{n}(W)^{+} = M_{n}(V)^{+}\cap M_{n}(W)_{sa}
\end{equation*}
then $(W, \{M_{n} (W)^{+}\})$ becomes a matrix ordered space and is
called a {\it matrix ordered subspace} of $(V, \{M_{n} (V)^{+}\})$.

A matrix ordered subspace $W$ of $V$ is called a {\it matrix order
ideal}, if for any $v\in V^{+}, v\in W^{+}$ whenever $v\leq w$ for some
$w\in W^{+}$. (In this case this property is carried over to
$M_{n}(W)^{+}$ for all $n$. 

Finally, let $(V, \{\left\|\right\|_{n}\}, \{M_{n}(V)^{+}\})$ be an
$L^{\infty}$-matricially Riesz normed space and suppose that $W$ is a
matrix order ideal in $V$ with an order unit $e$. If $(W, e)$ is a
matrix order unit space such that $e$ determines the matrix norm of $V$
on $W$ as the matrix order unit norm, then $(W, e)$ is called a matrix
order unit ideal of $V$.

\section{Order unit and matrix order ideal}

Let $(V, \{M_{n}(V)^{+}\})$ be a matrix ordered space and let $a\in V, a
\neq 0$. Put
\begin{equation*}
V_{a} =\left\lbrace 
v \in V: \left(\begin{array}{cc}
ka &v\\
v^{*} &ka
\end{array}\right)
\in M_{2}(V)^{+}, \textrm{for some}\ K>0
\right\rbrace.
\end{equation*}
Then $V_{a}$ becomes a matrix order ideal in $V$ and $a$ is an order
unit for $V_{a}$. Conversely, if $W$ is a matrix order ideal of $V$ and
$a\in W^{+}$ is an order unit for $W$, then $W = V_{a}$. In particular,
if $e$ is an order unit for $V$, then $V$ have no matrix order ideals
containing $e$, other than itself.

Now, let $(V, \{\left\|\right\|_{n}\}, \{M_{n}(V)^{+}\})$ be an
$L^{\infty}$-mRn space and suppose that $a\in V^{+}$ with
$\left\|a\right\|_{1} = 1$. We say that $a$ has {\it order unit
property} in $V$, if for each $v\in V_{a}$, we have
\begin{equation*}
\left(\begin{array}{cc}
\left\|v\right\|_{1}a &v\\
v^{*} &\left\|v\right\|_{1}a 
\end{array}
\right) \in M_{2}(V)^{+}.
\end{equation*}
If $a^{n}$ have order unit property for all $n$, we say that $a$ have
{\it matrix order unit property} in $V$. The following result is a
simple consequence of this definition.

\setcounter{defin}{0}
\begin{propo}\label{pro2.1}$\left.\right.$\vspace{.5pc}

\noindent Let $(V, \{\left\|\right\|_{n}\}, \{M_{n}(V)^{+}\})$ be an $L^{\infty}$-
mRn space and suppose that $a\in V^{+}$ with $\left\| v\right\|_{1} =
1$. Then $\hbox{`}a\hbox{'}$ have matrix order unit property in $V$ if
and only if $(V_{a}, a)$ is a matrix order unit ideal in $V$. In
particular, if $(W, e)$ is a matrix order unit ideal in $V${\rm ,} then
$\hbox{`}e\hbox{'}$ have matrix order unit property in $V$ and $W =
V_{e}$.
\end{propo}

\begin{rem}

{\rm It is possible to get an order unit $a\in V^{+}$ with $\left\|
a\right\|_{1} = 1$ which does not have order unit property in $V$. To
see this, consider the matrix order unit space $(M_{2}, I_{2})$. Then $a
= \left(\begin{array}{cc} 1 &0\\
0 &\displaystyle {1}/{2}
\end{array}\right) \in M_{2}^{+}$ and $\left\| a\right\|_{1} = 1$.
Moreover, $a$ determines the operator norm in the subspace defined as $A =
\{\alpha a: \alpha \in C\}$. Now, it is easy to see that $a$ is an order
unit for $M_{2}$. Finally, note that the operator norm of $I_{2}$ is 1
and its order unit norm determined by $a$ is 2. Thus $a$ does not have
order unit property in $M_{2}$.}
\end{rem}

\section{Order units in a $\hbox{\sa C}^{*}$-algebra}

Let $H$ be a complex Hilbert space and let $p\in B(H)$ be an
(orthogonal) projection. Define
\begin{align*}
B_{p}(H) & = \{x\in B(H): xp = px = x\}\\
& = \{pxp: x\in B(H)\}
\end{align*}
and
\begin{equation*}
B^{p}(H) = \{x\in B(H): x(H)\subset p(H)\quad \hbox{and}\quad
x^{*}(H)\subset p(H)\}.
\end{equation*}

\setcounter{defin}{0}
\begin{lemm}\label{lem3.1}
$B(H)_{p} = B_{p}(H) = B^{p}(H)\cong B(pH)$.\vspace{.2cm}
\end{lemm}

\begin{proof}
Evidently, $B_{p}(H) = B^{p}(H) \cong B(pH)$.

Let $x\in B(H)_{p}$. Then $\left(\begin{array}{cc}
\alpha p &x\\
x^{*} &\alpha p
\end{array}\right)\geq 0$ for some $\alpha >0$. Thus, for each pair
$\xi, \eta \in H$, we have
\begin{equation*}
\left\langle \left(\begin{array}{cc}
\alpha p &x\\
x^{*} &\alpha p
\end{array}
\right)\left(\begin{array}{c}
\xi\\
\eta
\end{array}\right)
\right\rangle\geq 0
\end{equation*}
or equivalently
\begin{equation*}
\alpha \langle p\xi, \xi\rangle + \langle x\eta, \xi \rangle + \langle
x^{*}\xi, \eta\rangle  + \alpha\langle \eta, \eta \rangle \geq 0.
\end{equation*}
In particular, for $\xi \in p(H)^{\bot}$ and $\eta \in H$, we have
\begin{equation*}
\langle x\eta, \xi \rangle + \langle x^{*}\xi, \eta \rangle + \alpha
\langle \eta, \eta \rangle\geq 0.
\end{equation*}
Fixing $\eta$ and replacing $\xi$ by $i^{k}\lambda \xi (k = 0, 1, 2, 3$
and $\lambda >0$), we conclude that $\langle x\eta, \xi \rangle\geq 0$
for all $\eta \in H$ and $\xi \in p(H)^{\bot}$. Thus $x(H)\subset p(H)$.
Next, interchanging roles of $\xi$ and $\eta$, we conclude that
$x^{*}(H)\subset p(H)$. Thus $B(H)_{p}\subset B^{p}(H)$.

Conversely, if $x\in B_{p}(H)$, then as $\left(\begin{array}{cc}
\left\| x\right\|I &x\\
x^{*} &\left\| x\right\|I
\end{array}\right) \geq 0$, we have
\begin{align*}
0 &\leq \left(\begin{array}{cc}
p &0\\
0 &p
\end{array}
\right) \left(\begin{array}{cc}
\left\|x \right\|I &x\\
x^{*} &\left\| x\right\|I 
\end{array}
\right) 
\left(\begin{array}{cc}
p &0\\
0 &p
\end{array}
\right)\\
\end{align*}
\begin{align*}
&=\left(\begin{array}{cc}
\left\| x\right\|p &pxp\\
px^{*}p &\left\|x\right\|p\end{array}
\right)
=\left(\begin{array}{cc}
\left\|x \right\|p &x\\
x^{*} &\left\|x\right\|p
\end{array}
\right)
\end{align*}
for $pxp = x$ and $x^{*} = (pxp)^{*} = px^{*}p$. Hence follows the
result.
\hfill $\ab$
\end{proof}

\setcounter{defin}{1}
\begin{cor}\label{cor3.2}$\left.\right.$\vspace{.5pc}

\noindent Let $A$ be a $C^{*}$-algebra and suppose that $p\in A^{+}$ is a
projection. Then $A_{p} = pAp = \{pxp: x\in A\}$. In particular{\rm ,} every
projection in $A$ has the matrix order unit property.
\end{cor}

\begin{theo}[\!]\label{the3.3}
Let $A$ be a $C^{*}$-algebra{\rm ,} $p\in A^{+}$ with $\left\| p\right\| = 1$.
If $p$ has the order unit property in $A${\rm ,} then $p$ is a projection.
\end{theo}

\begin{proof}
First we show that $A_{p}$ is a $C^{*}$-subalgebra of $A$. Let $x\in
A_{p}$. Then $\left(\begin{array}{cc}
kp &x\\
x^{*} &kp
\end{array}\right)\in M_{2}(A)^{+}$, for some $k>0$. Identifying $A$ in
a suitable $B(H)$, we conclude that $x^{*}x \leq k^{2}p$ (\cite{PV},
Exercise 3.2(ii)) so that $x^{*}x\in A_{p}$. It follows, from
polarization identity, that $A_{p}$ is a $^{*}$-algebra. To see that
$A_{p}$ is norm closed, let $\{x_{n}\}$ be any sequence in $A_{p}$
converging to $x\in A$. Let $\left\|x_{n}\right\|\leq k$ for all $n$ and
some $k>0$. Since $p$ has order unit property in $A$, we have $\left(\begin{array}{cc}
kp &x_{n}\\
x^{*}_{n} &kp
\end{array}\right)\in M_{2}(A)^{+}$ for all $n$. It follows, from 
norm-closedness of $M_{2}(A)^{+}$, that $\left(\begin{array}{cc}
kp &x\\
x^{*} &kp
\end{array}\right)\in M_{2}(A)^{+}$. Thus $A_{p}$ is a $C^{*}$-subalgebra of $A$.

Now identify $A_{p}$ in some $B(K)$ such that $A_{p}$ acts \hbox{non-degenerately}
on $K$. Let $\{u_{\lambda}\}$ be an approximate identity of
$A_{p}$. Then $u_{\lambda}\uparrow I$ where $I$ is the identity operator on
$K$ (see e.g.~\cite{PG}, 2.2.4). Now for each $\lambda,
\left\|u_{\lambda}\right\|\leq 1$ and $u_{\lambda}\geq 0$, so that
$u_{\lambda}\leq p$. Thus $I\leq p$. Again, as $\left\| p\right\| = 1$
and $p\geq 0$, we have $p\leq I$. Thus $p = I$. In particular, $A_{p}$
is a unital $C^{*}$-algebra. Hence $p$ is a projection in $A$.
\hfill $\ab$
\end{proof}

Now summing up Corollary~\ref{cor3.2} and Theorem~\ref{the3.3}, we
obtain the following characterization.

\begin{theo}[\!]\label{the3.4}

Let $A$ be a $C^{*}$-algebra and suppose that $p\in A^{+}$ with $\left\|
p\right\| = 1$. Then $p$ has the order unit property in $A$ if and only
if $p$ is a projection in $A$. In this case, $A_{p} = pAp$.
\end{theo}

\begin{rem}

{\rm (1) Matrix order unit ideals in a $C^{*}$-algebra are precisely unital
hereditary $C^{*}$-subalgebras. (2) Matrix order unit ideals in a 
von-Neumann algebra are precisely hereditary von-Neumann subalgebras.
(3) Matrix order unit ideals in $B(H)$ are precisely of the form $B(K)$
where $K$ is a closed subspace of $H$.}
\end{rem}


\begin{thebibliography}{99}

\bibitem{CM} Choi~M~D and Effros~E~G, Injectivity and operator spaces,
{\it J.~Funct.~Anal.} {\bf 24}~(1977) \hbox{156--209}

\bibitem{KA} Karn~A~K and Vasudevan~R, Matrix norms in matrix ordered
spaces, {\it Glasnik Mathemati\v{c}ki} {\bf 32(52)}~(1997) 87--97

\bibitem{KAK} Karn~A~K and Vasudevan~R, Approximate matrix order unit
spaces, {\it Yokohama Math.~J.} {\bf 44}~(1997) 73--91

\bibitem{KAK1} Karn~A~K and Vasudevan~R, Matrix duality for matrix
ordered spaces, {\it Yokohama Math.~J.}~{\bf 45}~(1998) 1--18

\bibitem{KAK2} Karn~A~K and Vasudevan~R, Characterizations of
matricially Riesz normed spaces, {\it Yokohama Math.~J.}~{\bf 47}~(2000) 143--153

\bibitem{PV} Paulsen~V~I, Completely bounded maps and dialations, {\it
Pitman Research Notes in Mathematics, Longman Scientific and Technical},
London, 1986 

\bibitem{PG} Pedersen~G~K, $C^{*}$-algebras and their group
automorphisms (Academic Press)~(1978)

\bibitem{RZ} Ruan~Z~J, Subspaces of $C^{*}$-algebras, {\it
J.~Funct.~Anal.}~{\bf 76}~(1988) 217--230 

\bibitem{SW} Schreiner~W~J, Matrix regular operator spaces, {\it
J.~Funct.~Anal.}~{\bf 152}~(1998) 136--175
\end{thebibliography}
\end{document}